\input amstex
\documentstyle{amsppt}
\magnification=\magstep1
 \hsize 13cm \vsize 18.35cm \pageno=1
\loadbold \loadmsam
    \loadmsbm
    \UseAMSsymbols
\topmatter
\NoRunningHeads
\title New approach to $q$-Euler, Genocchi numbers and their
interpolation functions
\endtitle
\author
  Taekyun Kim
\endauthor

\abstract In [1], Cangul-Ozden-Simsek constructed a $q$-Genocchi
numbers  of higher order  and gave Witt's formula of these numbers
by using  a $p$-adic fermionic integral on $\Bbb Z_p$.  In this
paper, we give another constructions of a $q$-Euler and Genocchi
numbers of higher order, which are different than their $q$-Genocchi
and Euler numbers of higher order. By using  our $q$-Euler and
Genocchi numbers of higher order, we can investigate  the
interesting relationship between $q$-$w$-Euler numbers and
$q$-$w$-Genocchi numbers. Finally, we give the interpolation
functions of these numbers.

\endabstract
\thanks  2000 AMS Subject Classification: 11B68, 11S80
integral
\newline  The present Research has been conducted by the research
Grant of Kwangwoon University in 2008
\endthanks
\endtopmatter

\document

{\bf\centerline {\S 1. Introduction/ Preliminaries}}

 \vskip 20pt
In [1], Cangul-Ozden-Simsek presented a systematic study of some
families of the $q$-Genocchi numbers of higher order. By applying
the $q$-Volkenborn integral in the sense of fermionic, they also
constructed the $q$-extension of Euler zeta function, which
interpolates these numbers at negative integers. Some properties of
Genocchi numbers of higher order were treated in [1, 11, 12, 13, 14,
15, 16, 17]. So, we will construct the appropriate the $q$-Genocchi
numbers of higher order to be related $q$-Euler numbers of higher
order for doing  study the $q$-extension of Genocchi numbers of
higher order in this paper. Let $p$ be a fixed odd prime. Throughout
this paper $\Bbb Z_p$, $\Bbb Q_p$, $\Bbb C$ and $\Bbb C_p$ will,
respectively, denote the ring of $p$-adic rational integers, the
field of $p$-adic rational numbers, the complex number field and the
completion of algebraic closure of $\Bbb Q_p$. Let $v_p$ be the
normalized exponential valuation of $\Bbb C_p$ with
$|p|_p=p^{-v_p(p)}=\frac{1}{p}.$ The $q$-basic natiral numbers are
defined by $[n]_q=\frac{1-q^n}{1-q}=1+q+q^2+\cdots+q^{n-1}$ for $n
\in\Bbb N$, and the binomial coefficient is defined as
$$\binom{n}{k}=\frac{n!}{(n-k)!k!}=\frac{n(n-1)\cdots(n-k+1)}{k!}.$$
The binomial formulas are well known that
$$(1-b)^n=\sum_{i=0}^n \binom{n}{i}(-1)^i b^i, \text{ and }
\frac{1}{(1-b)^n}=\sum_{i=0}^n\binom{n+i-1}{i}b^i, \text{ (see
[1-20])}.$$ In this paper, we use the notation
$$[x]_q=\frac{1-q^x}{1-q}, \text{ and }
[x]_{-q}=\frac{1-(-q)^x}{1+q}, \text{ (see [1-22])}.$$

We say that $f$ is uniformly differentiable function at a point
$a\in\Bbb Z_p$, and write $f\in UD(\Bbb Z_p),$ if the difference
quotient $F_{f}(x, y)=\frac{f(x)-f(y)}{x-y}$ have a limit
$f^{\prime}(a)$ as $(x,y)\rightarrow (a,a).$ For $f \in UD(\Bbb
Z_p),$ $q\in \Bbb C_p$ with $|1-q|_p<1 ,$ an invariant $p$-adic
$q$-integral is defined as
$$I_{-q}(f)=\int_{\Bbb Z_p}f(x) d\mu_{-q}(x)=\lim_{N\rightarrow
\infty}\frac{1+q}{1+q^{p^N}}\sum_{x=0}^{p^N-1}f(x)(-q)^x, \text{ see
[9]}. \tag1$$ Thus, we have the following integral relation:
$$qI_{-q}(f_1)+I_{-q}(f)=(1+q)f(0), \text{ where $f_{1}(x)=f(x+1)$ }.$$
The fermionic $p$-adic  invariant integral on $\Bbb Z_p$ is defined
as
$$I_{-1}(f)=\lim_{q\rightarrow 1}I_{-q}(f)=\int_{\Bbb Z_p}f(x)d\mu_{-1}(x).$$
For $n\in \Bbb N$, let $f-n(x)=f(x+n)$. From (1), we can derive
$$q^nI_{-q}(f_n)=(-1)^nI_{-q}(f)+[2]_q\sum_{l=0}^{n-1}(-1)^{n-1-l}q^lf(l).$$
It is known that the $w$-Euler polynomials are defined as
$$\frac{2e^{xt}}{we^t+1}=\sum_{n=0}^{\infty}E_{n,w}(x)\frac{t^n}{n!},
\text{ (see [ 20, 21, 22])}. \tag2$$ Note that $E_{n,w}(0)=E_{n,w}$
are called the $w$-Euler numbers. $w$-Genocchi polynomials are
defined as
$$\frac{2t}{we^t+1}e^{xt}=\sum_{n=0}^{\infty}G_{n,w}(x)\frac{t^n}{n!},
\text{ (see [1])}.$$ In the special case $x=0$, $G_{n,w}(0)=G_{n,w}$
are called $w$-Genocchi numbers. In [1], Cangul-Ozden-Simsek have
studied $w$-Genocchi numbers of order $r$ as follows. For $w\in \Bbb
C_p$ with $|1-w|_p<1,$ the $w$-Genocchi numbers of order $r$ are
given by
$$\int_{\Bbb Z_p}\cdots \int_{\Bbb
Z_p}e^{t(x_1+\cdots+x_r)}d\mu_{-1}(x_1)\cdots d\mu_{-1}(x_r)
=2^r\left(\frac{t}{we^t+1}\right)^r=\sum_{n=0}^{\infty}G_{n,w}^{(r)}\frac{t^n}{n!}.\tag3$$
They also considered  the $q$-extension of (3) as follows.

$$\aligned
&\int_{\Bbb Z_p}\cdots\int_{\Bbb
Z_p}q^{\sum_{v=0}^r(h-v)x_v}e^{t(x_1+\cdots+x_r)}d\mu_{-1}(x_1)\cdots
d\mu_{-1}(x_r)\\
&=\frac{2^rt^r}{(q^he^t+1)(q^{h-1}e^t+1)\cdots(q^{h-r+1}e^t+1)}=\sum_{n=0}^{\infty}G_{n,q}^{(h,r)}\frac{t^n}{n!},
\text{ (see [1]) }.
\endaligned\tag4$$
 The purpose of this paper is to give
another construction of $q$-Euler numbers and $q$-Genocchi numbers
of higher order, which are different than a  $q$-Genocchi numbers of
Cangul-Ozden-Simsek.

\vskip 20pt

{\bf\centerline {\S 2.New approach to $q$-Euler, Genocchi numbers
and polynomials}} \vskip 10pt

In this section, we assume that $q\in \Bbb C_p$ with $|1-q|_p<1$ and
let $w\in\Bbb C_p$ with $|1-w|_p<1.$ The $w$-Euler polynomials of
order $r$, denoted $E_n^{(r)}(x)$, are defined as
$$e^{xt}\left(\frac{2}{we^t+1}\right)^r=
\underbrace{\left(\frac{2}{we^t+1}\right)\times\cdots\times\left(\frac{2}{we^t+1}\right)}_{r-times}e^{xt}
=\sum_{n=0}^{\infty} E_{n,w}^{(r)}(x)\frac{t^n}{n!}, \text{ (see
[1])}.$$ The values of $E_{n,w}^{(r)}(x)$ at $x=0$ are called
$w$-Euler number of order $r$: when $r=1$ and $w=1$, the polynomials
or numbers are called the ordinary Euler polynomials or numbers.
When $x=0$ or $r=1$, we use the following notation: $ E_{n,w}^{(r)}$
denote $E_{n,w}^{(r)}(0)$, $E_{n,w}(x)$ denote $E_{n,w}^{(1)}(x),$
and $E_{n,w}$ denote $E_{n,w}^{(1)}(0)$.

 It is known that
$$\int_{\Bbb Z_p}w^xe^{tx}
d\mu_{-1}(x)=\frac{2}{we^t+1}=\sum_{n=0}^{\infty}E_{n,w}\frac{t^n}{n!},
$$
and
$$\int_{\Bbb Z_p}w^ye^{t(x+y)}
d\mu_{-1}(y)=\frac{2}{we^t+1}e^{xt}=\sum_{n=0}^{\infty}E_{n,w}(x)\frac{t^n}{n!},
\text{ (see ([1])}
$$
The higher order $w$-Euler numbers and polynomials are given by
$$\int_{\Bbb Z_p}\cdots\int_{\Bbb
Z_p}w^{x_1+\cdots+x_r}e^{(x_1+\cdots+x_r)t}d\mu_{-1}(x_1)\cdots
d\mu_{-1}(x_r)=\left(\frac{2}{we^{t}+1}
\right)^r=\sum_{n=0}^{\infty}E_{n,w}^{(r)}\frac{t^n}{n!}, $$ and
$$
\int_{\Bbb Z_p}\cdots\int_{\Bbb Z_p}w^{\sum_{i=1}^r
x_i}e^{(\sum_{i=1}^r x_i+x)t}\prod_{i=1}^r d\mu_{-1}(x_i)
=\left(\frac{2}{we^t+1} \right)^re^{xt}
=\sum_{n=0}^{\infty}E_{n,w}(x)^{(r)}\frac{t^n}{n!}. \tag5$$

Thus, we have the following theorem.

 \proclaim{ Proposition 1}
For $n\in\Bbb N$, we have
$$E_{n,w}^{(r)}(x)=\int_{\Bbb Z_p}\cdots \int_{\Bbb Z_p}w^{x_1+\cdots+x_r}(x_1 +\cdots+x_r+x)^nd\mu_{-1}(x_1)\cdots d\mu_{-1}(x_r),\tag6$$
and
$$E_{n,w}^{(r)}=\int_{\Bbb Z_p}\cdots \int_{\Bbb Z_p}w^{x_1+\cdots+x_r}(x_1 +\cdots+x_r)^nd\mu_{-1}(x_1)\cdots d\mu_{-1}(x_r).$$
\endproclaim
From the results of Cangul-Ozden-Simsek (see [1]), we can derive the
following $w$-Genocchi polynomials of order $r$ as follows.

$$\aligned
& t^r\int_{\Bbb Z_p}\cdots \int_{\Bbb Z_p}
w^{x_1+x_2+\cdots+x_r}e^{(x_1+x_2+\cdots+x_r+x)t}d\mu_{-1}(x_1)d\mu_{-1}(x_2)\cdots
d\mu_{-1}(x_r)\\
&=\left( \frac{2}{we^t
+1}\right)^re^{xt}=\sum_{n=0}^{\infty}G_{n,w}^{(r)}(x)\frac{t^n}{n!},
\text{ (see [1]).}
\endaligned\tag7$$
From (7), we note that
$$\aligned
& t^r\int_{\Bbb Z_p}\cdots \int_{\Bbb Z_p}
w^{x_1+x_2+\cdots+x_r}e^{(x_1+x_2+\cdots+x_r+x)t}d\mu_{-1}(x_1)d\mu_{-1}(x_2)\cdots
d\mu_{-1}(x_r)\\
&=\sum_{n=0}^{\infty}\int_{\Bbb Z_p}\cdots\int_{\Bbb
Z_p}w^{x_1+\cdots+x_r}(x_1+\cdots+x_r+x)^n
r!\binom{n+r}{r}\frac{t^{n+r}} {(n+r)!}.
\endaligned\tag8$$

By (6), (7) and (8), we easily see that
$G_{0,w}^{(r)}(x)=G_{1,w}^{(r)}(x)=\cdots=G_{n-1,w}^{(r)}(x)=0,$ and
$$\aligned
\frac{G_{n+r,w}^{(r)}(x)}{r!\binom{n+r}{r}}&=\int_{\Bbb Z_p}\cdots
\int_{\Bbb Z_{p}}w^{x_1+\cdots+x_r}(x_1+x_2+\cdots+x_r+x)^n
d\mu_{-1}(x_1)\cdots d\mu_{-1}(x_r)\\
&=E_{n,w}^{(r)}(x).
\endaligned$$

In the viewpoint of the $q$-extension of (6), let us define the
$w$-$q$-Euler numbers of order $r$ as follows.
$$E_{n,w,q}^{(r)}=\int_{\Bbb Z_p}\cdots\int_{\Bbb
Z_p}w^{x_1+\cdots+x_r}[x_1+\cdots+x_r]_q^nq^{-(x_1+\cdots+x_r)}d\mu_{-q}(x_1)\cdots
d\mu_{-q}(x_r).\tag9$$

From (9), we note that
$$\aligned
E_{n,w,q}^{(r)}&=\frac{1}{(1-q)^n}\sum_{l=0}^n\binom{n}{l}(-1)^l\left(\frac{[2]_q}{1+wq^l}\right)^r\\
&=[2]_q^r\frac{1}{(1-q)^n}\sum_{l=0}^n
\binom{n}{l}(-1)^l\sum_{m=0}^{\infty}\binom{m+r-1}{m}(-1)^mw^mq^{lm}\\
&=[2]_q^r\sum_{m=0}^{\infty}\binom{m+r-1}{m}(-1)^mw^m[m]_q^n.
\endaligned\tag10$$
Therefore, we obtain the following theorem.

\proclaim{ Theorem 2} For $n\in\Bbb N$, we have

$$E_{n, w,
q}^{(r)}=[2]_q^r\sum_{m=0}^{\infty}\binom{m+r-1}{m}(-1)^mw^m[m]_q^n.$$
\endproclaim
Let $F(t,x|q)=\sum_{n=0}^{\infty}E_{n,w,q}^{(r)}\frac{t^n}{n!}.$
Then we have
$$\aligned
F(t,w|q)&=\int_{\Bbb Z_p}\cdots \int_{\Bbb Z_p}w^{x_1+\cdots+x_r}
q^{-(x_1+\cdots+x_r)}e^{[x_1+\cdots+x_r]_qt}d\mu_{-q}(x_1)\cdots
d\mu_{-q}(x_r)\\
&=[2]_q^r\sum_{m=0}^{\infty}\binom{m+r-1}{m}(-1)^m w^me^{[m]_qt}.
\endaligned$$

Thus, we obtain the following corollary.

\proclaim{ Corollary 3} Let
$F(t,x|q)=\sum_{n=0}^{\infty}E_{n,w,q}^{(r)}\frac{t^n}{n!}.$ Then we
have
$$F(t,x|q)=[2]_q^r\sum_{m=0}^{\infty}\binom{m+r-1}{m}(-1)^m
w^me^{[m]_qt}.$$
\endproclaim

Note that $$\frac{d^k
}{{dt}^k}F(t,x|q)\big|_{t=0}=E_{k,w,q}^{(r)}=[2]_q^r\sum_{m=0}^{\infty}\binom{m+r-1}{m}(-1)^mw^m[m]_q^k.$$

Let us define the $q$-extension of $w$-Genocchi numbers of order $r$
as follows.
$$\aligned
&t^r\int_{\Bbb Z_p}\cdots\int_{\Bbb Z_p}q^{-(
x_1+\cdots+x_r)}w^{x_1+\cdots+x_r}e^{[x_1+\cdots+x_r]_qt}d\mu_{-q}(x_1)\cdots
d\mu_{-q}(x_r)\\
&=t^r[2]_q^r\sum_{m=0}^{\infty}\binom{m+r-1}{m}(-1)^m w^m
e^{[m]_qt}=\sum_{n=0}^{\infty}G_{n,w,q}^{(r)}\frac{t^n}{n!}.
\endaligned\tag11$$
From (11), we can derive
$$\aligned
&\sum_{n=0}^{\infty}\int_{\Bbb Z_p}\cdots\int_{\Bbb Z_p}
q^{-(x_1+\cdots+x_r)}w^{x_1+\cdots+x_r}[x_1+\cdots+x_r]_q^n
\prod_{i=1}^r d\mu_{-q}(x_i)\frac{r!\binom{n+r}{r}t^{n+r}}{(n+r)!}\\
&=\sum_{n=0}^{\infty}G_{n,w,q}^{(r)}\frac{t^n}{n!}=\sum_{n=0}^{\infty}
G_{n+r, w,q}^{(r)}\frac{t^{n+r}}{(n+r)!}.
\endaligned\tag12$$
By (9), (11) and (12), we obtain the following theorem.

 \proclaim{Theorem 4}
For $n \in \Bbb Z_{+}$, and $r\in\Bbb N$,  we have
$$\aligned
\frac{G_{n+r,w,q}^{(r)}}{r!\binom{n+r}{r}}&=\int_{\Bbb
Z_p}\cdots\int_{\Bbb
Z_p}q^{-(x_1+\cdots+x_r)}w^{x_1+\cdots+x_r}[x_1+\cdots+x_r]_q^n
d\mu_{-q}(x_1)\cdots d\mu_{-q}(x_r)\\
 &=E_{n,w,q}^{(r)},
\endaligned$$
and $G_{0,w,q}^{(r)}=G_{1,w,q}^{(r)}=\cdots=G_{r-1,w,q}^{(r)}=0.$
 \endproclaim

In the sense of the q-extension of (6), we can consider the
$w$-$q$-Euler polynomials of order $r$ as follows.
$$\aligned
E_{n,w,q}^{(r)}(x)&=\int_{\Bbb Z_p}\cdots\int_{\Bbb
Z_p}w^{x_1+\cdots+x_r}q^{-(x_1+\cdots+x_r)}[x_1+\cdots+x_r]_q^n
d\mu_{-q}(x_1)\cdots d\mu_{-q}(x_r)\\
&=\frac{1}{(1-q)^n}\sum_{l=0}^n\binom{n}{l}(-1)^lq^{lx}\left(\frac{[2]_q}{1+wq^l}\right)^r\\
& =[2]_q^r\sum_{m=0}^{\infty}\binom{m+r-1}{m}(-1)^m[m+x]_q^nw^m.
\endaligned\tag13$$

Let $F(t,w,
x|q)=\sum_{n=0}^{\infty}E_{n,w,q}^{(r)}(x)\frac{t^n}{n!}. $ Then we
have
$$\aligned
&F(t,
w,x|q)=[2]_q^r\sum_{m=0}^{\infty}\binom{m+r-1}{m}(-1)^mw^me^{[m+x]_qt}\\
&=\int_{\Bbb Z_p}\cdots\int_{\Bbb
Z_{p}}w^{x_1+\cdots+x_r}q^{-(x_1+\cdots+x_r)}e^{[x_1+\cdots+x_r+x]_qt}d\mu_{-q}(x_1)\cdots
d\mu_{-q}(x_r).
\endaligned\tag14$$

In the viewpoint of (7), we can define the $w$-$q$-Genocchi
polynomials of order $r$ as follows.

$$\sum_{n=0}^{\infty}G_{n,w,q}^{(r)}(x)\frac{t^n}{n!}=[2]^rt^r\sum_{m=0}^{\infty}\binom{m+r-1}{m}(-1)^mw^me^{[m+x]_qt}.
\tag15$$

By (14) and (15), we obtain the following theorem.

\proclaim{Theorem 5} For $n\in\Bbb Z_{+}$, $r\in \Bbb N$, we have
$$\aligned
&\frac{G_{n+r,w,q}^{(r)}(x)}{r!\binom{n+r}{r}}=E_{n,w,q}^{(r)}(x)\\
&=\int_{\Bbb Z_p}\cdots\int_{\Bbb
Z_p}w^{x_1+\cdots+x_r}q^{-(x_1+\cdots+x_r)}[x_1+\cdots+x_r]_q^n
d\mu_{-q}(x_1)\cdots d\mu_{-q}(x_r),
\endaligned$$
and
$$G_{0,w,q}^{(r)}(x)=\cdots=G_{r-1,w,q}^{(r)}(x)=0.$$
\endproclaim

\vskip 20pt

{\bf\centerline {\S 3. Further Remarks and Observation for the
multiple $w$-$q$-zeta function}} \vskip 10pt

In this section, we assume that $q\in\Bbb C$ with $|q|<1.$ For
$s\in\Bbb C$, $w=e^{2\pi i\xi}$($\xi\in \Bbb R$), let us define  the
Lerch type $q$-zeta function of order $r$ as follows.

$$\zeta_{q,w}^{(r)}(s)=[2]_q^r\sum_{m=1}^{\infty}\frac{\binom{m+r-1}{m}(-1)^mw^m}{[m]_q^s}.
\tag16$$ For $k\in \Bbb N$, by (10) and (16), we have
$\zeta_{q,w}^{(r)}(-k)=E_{k,w,q}^{(r)}.$

We now consider the $q$-extension of Hurwitz-Lerch type zeta
function of order $r$ as follows.
$$\zeta_{q,w}^{(r)}(s,x)=[2]_q^r\sum_{m=0}^{\infty}\frac{\binom{m+r-1}{m}(-1)^mw^m}{[m+x]_q^s},$$
where $x\in\Bbb C$ with $x\neq 0, -1, -2, \cdots,$ and $s\in\Bbb C$.

From (13) and (17), we can also derive
$$\zeta_{q,w}^{(r)}(-k, x)=E_{n, w,q}^{(r)}(x), \text{  for
$k\in\Bbb N$}.$$

\vskip 20pt

\quad Taekyun Kim

\quad Division of General Education-Mathematics, Kwangwoon
University, Seoul 139-701, S. Korea \quad e-mail:\text{
tkkim$\@$kw.ac.kr}


\Refs \widestnumber\key{999999}

\ref \key 1
 \by I. N. Cangul, H. Ozden, Y. Simsek
  \paper A new Approach to $q$-Genocchi numbers and their
 interpolations
 \jour  Nonlinear Analysis(2008),
 \yr
\pages \vol doi:10.1016/j.na.2008.11.040 \endref


\ref \key 2
 \by   L. Comtet
 \book Advanced combinatories, Reidel
 \publ  Dordrecht
 \yr 1974
 \endref

\ref \key 3
 \by  E.Deeba, D.Rodriguez
 \paper   Stirling's  series and Bernoulli numbers
 \jour  Amer. Math. Monthly
 \yr 1991
\pages 423-426 \vol 98 \endref

\ref \key 4
 \by  M. Cenkci, M. Can and V. Kurt
  \paper  $p$-adic interpolation functions and Kummer-type congruences for $q$-twisted
  Euler numbers
 \jour  Adv. Stud. Contemp. Math.
\yr 2004 \pages 203--216 \vol 9 \endref

\ref \key 5
 \by  F. T. Howard
  \paper Application of a recurrence for the Bernoulli numbers
 \jour  J. Number Theory
\yr 1995 \pages 157-172 \vol 52\endref

\ref \key 6
 \by  T. Kim
  \paper The modified $q$-Euler numbers and polynomials
    \jour Adv. Stud. Contemp. Math.
\yr 2008 \pages 161-170 \vol 16 \endref

\ref \key 7
 \by  T. Kim
  \paper Euler numbers and polynomials associated with zeta functions
  \jour  Abstract and Applied Analysis
\yr 2008 \pages 11 pages(Article ID 581582 )  \vol 2008 \endref

\ref \key 8
 \by  K. Shiratani, S. Yamamoto
  \paper  On a $p$-adic interpolation function for the Euler numbers
  and its derivatives
 \jour  Mem. Fac. Sci., Kyushu University Ser. A
 \yr 1985
\pages 113-125\vol 39\endref

\ref \key 9
 \by   H.J.H. Tuenter
  \paper A Symmetry of power sum polynomials and Bernoulli numbers
   \jour Amer. Math. Monthly
 \yr 2001
\pages 258-261\vol 108 \endref

\ref \key 10
 \by  T. Kim
  \paper  $q-$Volkenborn integration
 \jour  Russ. J. Math. Phys.
 \yr 2002
\pages 288--299 \vol 9 \endref

\ref \key 11
 \by  T. Kim
  \paper A Note on $p$-Adic $q$-integral on  $\Bbb Z_p$
  Associated with $q$-Euler Numbers
 \jour Adv. Stud. Contemp. Math.
 \yr 2007
\pages 133--138 \vol 15 \endref

\ref \key 12
 \by  T. Kim
  \paper  On $p$-adic interpolating function for $q$-Euler numbers and
   its derivatives
 \jour J. Math. Anal. Appl.
 \yr 2008
\pages  598--608\vol 339 \endref

\ref \key 13
 \by T. Kim
  \paper  $q$-Extension of the Euler formula and trigonometric functions
 \jour    Russ. J. Math. Phys.\yr 2007
\pages  275--278 \vol 14\endref

\ref \key 14
 \by T. Kim
  \paper  Power series and asymptotic series
  associated with the $q$-analog of the two-variable
  $p$-adic $L$-function
 \jour    Russ. J. Math. Phys.\yr 2005
\pages  186--196 \vol 12\endref

\ref \key 15
 \by T. Kim
  \paper  Non-Archimedean $q$-integrals associated
  with multiple Changhee $q$-Bernoulli polynomials
 \jour    Russ. J. Math. Phys.\yr 2003
\pages 91--98 \vol 10\endref

\ref \key 16
 \by   T. Kim
  \paper  $q$-Euler numbers and polynomials associated with $p$-adic $q$-integrals
 \jour  J. Nonlinear Math. Phys.
 \yr 2007
\pages 15--27 \vol 14 \endref


\ref \key 17
 \by   H. Ozden, Y. Simsek, S.-H. Rim, I.N. Cangul
  \paper  A note on $p$-adic $q$-Euler measure
 \jour  Adv. Stud. Contemp. Math.
 \yr 2007
\pages 233--239 \vol 14 \endref

\ref \key 18
 \by M. Schork,
  \paper Ward's "calculus of sequences", $q$-calculus and the limit $q\to-1$
   \jour Adv. Stud. Contemp. Math.
 \yr 2006
\pages 131--141 \vol 13 \endref

\ref \key 19 \by M. Schork
  \paper Combinatorial aspects of normal
  ordering and its connection to $q$-calculus
   \jour Adv. Stud. Contemp. Math.
 \yr 2007
\pages 49-57 \vol 15 \endref

\ref \key 20
 \by  Y. Simsek
  \paper On $p$-adic twisted $q\text{-}L$-functions related to generalized twisted Bernoulli numbers
   \jour  Russ. J. Math. Phys.
 \yr 2006
\pages 340--348 \vol 13 \endref

\ref \key 21
 \by  Y. Simsek
  \paper  Theorems on twisted $L$-function and twisted Bernoulli
  numbers
   \jour Advan. Stud. Contemp. Math.
 \yr 2005
\pages 205--218 \vol 11 \endref

\ref \key 22
 \by   Y. Simsek
  \paper $q$-Dedekind type sums related to $q$-zeta function and basic $L$-series
   \jour J. Math. Anal. Appl.
 \yr 2006
\pages 333-351\vol 318 \endref



\endRefs
\enddocument